\documentclass[11pt]{article}
\usepackage{amsmath,amssymb,amsthm,bm,graphicx,color,epsfig}

\renewcommand{\i}{\imath}

\begin{document}

\title{Fast Computation of Partial Fourier Transforms}
\author{Lexing Ying$^{\dagger}$ and Sergey Fomel$^{\sharp}$ \\
  \vspace{0.1in}\\
  $\dagger$ Department of Mathematics, University of Texas, Austin, TX 78712\\
  $\sharp$ Bureau of Economic Geology, University of Texas, Austin, TX 78712}

\date{January 2008}
\maketitle

\begin{abstract}
  We introduce two efficient algorithms for computing the partial
  Fourier transforms in one and two dimensions. Our study is motivated
  by the wave extrapolation procedure in reflection seismology. In
  both algorithms, the main idea is to decompose the summation domain
  of into simpler components in a multiscale way. Existing fast
  algorithms are then applied to each component to obtain optimal
  complexity. The algorithm in 1D is exact and takes $O(N\log^2 N)$
  steps. Our solution in 2D is an approximate but accurate algorithm
  that takes $O(N^2 \log^2 N)$ steps. In both cases, the complexities
  are almost linear in terms of the degree of freedom. We provide
  numerical results on several test examples.
\end{abstract}



\section{Introduction}
\label{sec:intro}

In this paper, we introduce efficient algorithms for the following
partial Fourier transform problem in one and two dimensions. In 1D,
let $N$ be a large integer and $c^0(t)$ be a smooth function on
$[0,1]$ with $ 0 \le c^0(t) \le 1$. We define a sequence of integers
$\{c_x, 0\le x < N\}$ by $c_x = N \lceil c^0(x/N) \rceil$. Given a
sequence of $N$ numbers $\{f_k, 0\le k<N\}$, the problem is to compute
$\{u_x, 0\le x < N\}$ where
\begin{equation}
u_x = \sum_{k <c_x} e^{2\pi\i xk/N} f_k,
\label{eq:pft1d}
\end{equation}
where we use $\i$ for $\sqrt{-1}$ throughout this paper. In 2D, the
problem is defined in a similar way. Now let $c^0(t)$ be a smooth
function on $[0,1]^2$ with $0 \le c^0(t) \le 1$. The array $\{c_x,
0\le x_1,x_2<N\}$ where $x=(x_1,x_2)$ is defined by $c_x = N \lceil
c^0(x/N) \rceil$. Given an array of $N^2$ numbers $\{f_k, 0\le
k_1,k_2<N\}$ where $k=(k_1,k_2)$, we want to compute $\{u_x, 0\le
x_1,x_2<N\}$ defined by
\begin{equation}
u_x = \sum_{|k|<c_x} e^{2\pi\i x\cdot k/N} f_k.
\label{eq:pft2d}
\end{equation}

Due to the existence of the summation constraints on $k$ in
\eqref{eq:pft1d} and \eqref{eq:pft2d}, the fast Fourier transform
cannot be used directly here. Direct computation of \eqref{eq:pft1d}
and \eqref{eq:pft2d} has quadratic complexity; i.e., $O(N^2)$
operations for \eqref{eq:pft1d} and $O(N^4)$ for \eqref{eq:pft2d}.
This can be expensive for large values of $N$, especially in 2D. In
this paper, we propose efficient solutions which have almost linear
complexity. Our algorithm for the 1D case is exact and takes
$O(N\log^2 N)$ steps, while our solution to the 2D problem is an
accurate approximate algorithm with an $O(N^2 \log^2 N)$ complexity.
We define the {\em summation domain} $D$ to be the set of all pairs
$(x,k)$ appeared in the summation, namely $D = \{ (x,k), k <c_x\}$ in
1D and $D = \{ (x,k), |k|<c_x\}$ in 2D. The main idea behind both
algorithms is to partition the summation domain into simple components
in a multiscale fashion. Fast algorithms are then invoked on each
simple component to achieve maximum efficiency.

The partial Fourier transform appears naturally in several settings.
The one which motivates our research is the one-way wave extrapolation method in
seismology \cite{biondi-2006-3si}, where one often needs to
compute an approximation to the following integral \cite{margrave-1999-wenps}:
\begin{equation}
\label{eq:seismic}
u_z(x) = \int e^{2\pi\i (x\cdot k + \sqrt{ \omega^2/v^2(x) - k^2 }\cdot
  z )} \widehat{u_0}(k) d k,
\end{equation}
where $\omega$ and $z$ are fixed constants (frequency and extrapolation depth), $v(x)$ is a given function (layer velocity),
and $\widehat{u_0}(k)$ is the Fourier transform of a function
$u_0(x)$. The wave modes that correspond to $|k| \le \omega/v(x)$ are
propagating waves, while the ones that correspond to $|k| \ge
w/\rho(x)$ are called evanescent. For the purposes of seismic imaging, one is often
interested in only the propagating mode and, therefore, we have the
following restricted integral to evaluate:
\[
\int_{|k|\le \omega/v(x)} e^{2\pi\i (x \cdot k + \sqrt{
    \omega^2/v^2(x) - k^2 }\cdot z )} \widehat{u_0}(k) d k.
\]
To make the computation efficient, the term with the square root under is often approximated with a functional form
\begin{equation}
\label{eq:approx}
e^{2\pi\i \sqrt{
    \omega^2/v^2(x) - k^2 }\cdot z} \approx \sum_n f_n(x)\,\psi_n(k)
\end{equation}
with a limited number of terms. The integral then reduces to
\[
\sum_n f_n(x)\,\left[
\int_{|k|\le \omega/v(x)} e^{2\pi\i x} \widehat{w_n}(k) d k.
\right]
\]
where $\widehat{w_n}(k) = \psi_n(k)\,\widehat{u_0}(k)$. The kernel of this approximation involves computing
\[
\int_{|k|\le w/\rho(x)} e^{2\pi\i x\cdot k} \widehat{w_n}(k) d k,
\]
which takes the forms of \eqref{eq:pft1d} and \eqref{eq:pft2d} after
discretization.

The rest of this paper is organized as follows. In Section
\ref{sec:1d}, we describe our 1D algorithm. The 2D algorithm is
presented in Section \ref{sec:2d}. Finally, we conclude with some
discussions for future work in Section \ref{sec:conc}. Throughout this
paper, we assume for simplicity that $N$ is an integer power of $2$.

\section{Partial Fourier Transform in 1D}
\label{sec:1d}

\subsection{Algorithm description}
\begin{figure}[ht]
  \begin{center}
    \includegraphics[height=2in]{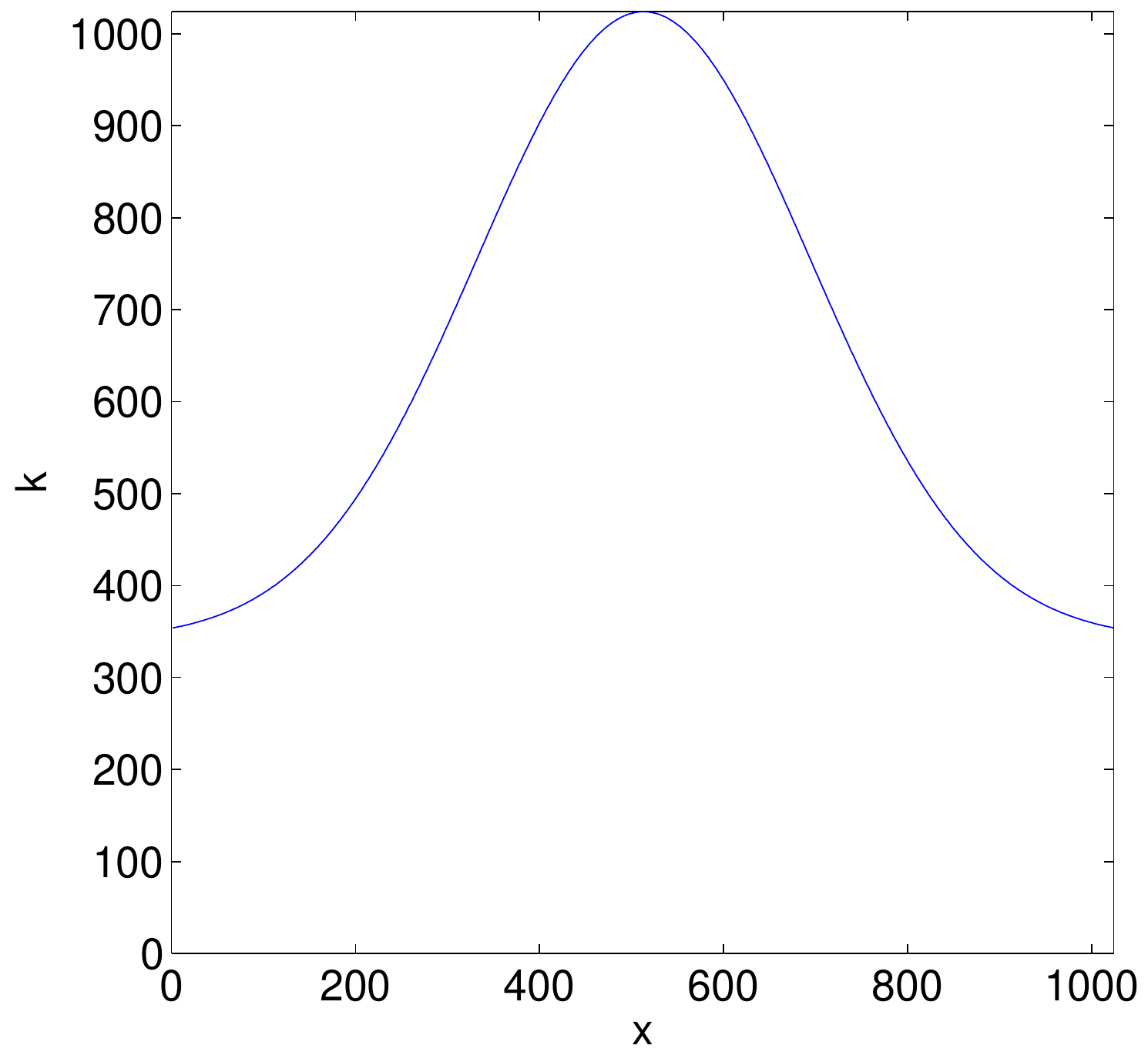}
    \includegraphics[height=2in]{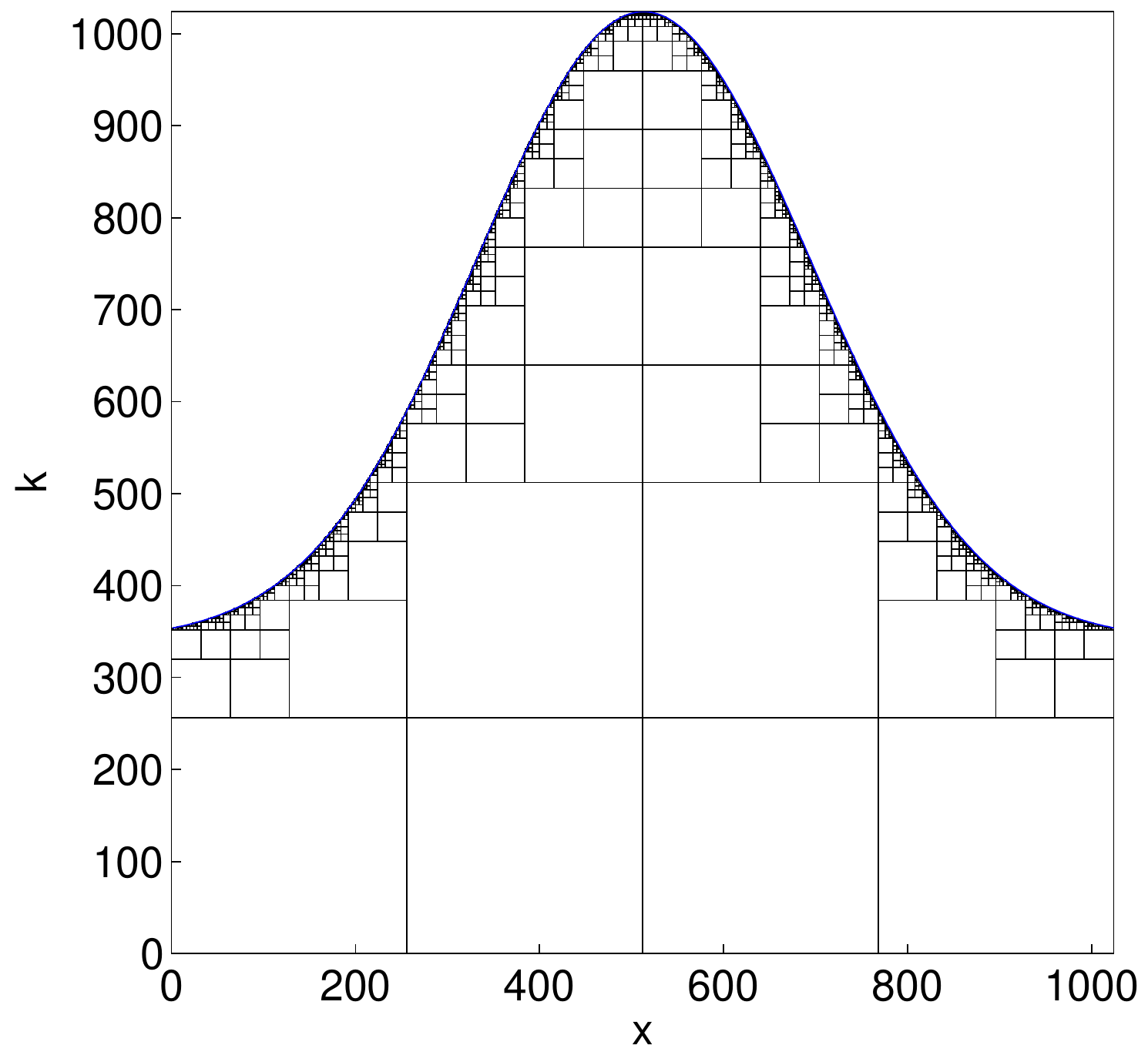}
  \end{center}
  \caption{Left: the curve stands for $\{c_x\}$ and the summation
    domain $D$ is the region below the curve. Right: the multiscale
    decomposition constructed for $D$. The boxes increase their size
    geometrically as they move away from the curve.}
  \label{fig:decomp1d}
\end{figure}

We define the summation domain by $D = \{(x,k), k<c_x\}$. The idea
behind our algorithm for computing
\[
u_x = \sum_{k <c_x} e^{2\pi\i xk/N} f_k,\quad 0\le x < N
\]
is to decompose $D$ in a multiscale fashion. Starting from the top
level box $[0,N]^2$, we partition the domain recursively. If a box $B$
is fully inside $D$, it is not subdivided and we keep it inside the
decomposition. If a box $B$ is fully outside $D$, it is discarded.
Finally, if a box $B$ has parts that belong to both $D$ and $[0,N]^2
\setminus D$, it is further subdivided into four child boxes with
equal size. At the end of this process, our decomposition contains a
group of boxes with dyadic sizes and the union of these boxes is
exactly equal to $D$ (see Figure \ref{fig:decomp1d}).

Let us consider a single box $B$ of the decomposition. Suppose that
$B$ is of size $s$ and that its lower-left corner of $B$ is
$(x^B,k^B)$. The part of the summation that associated with $B$ is
\begin{equation}
  \sum_{k^B\le k < k^B+s} e^{2\pi\i xk/N} f_k
  \label{eq:compB}
\end{equation}
for each $x^B \le x <x^B+s$. Denoting $x=x^B+x'$ and $k=k^B+k'$, we
can write this into a matrix form $M f$ with
\[
(M)_{x'k'} = e^{2\pi\i (x^B+x')(k^B+k')/N} =
e^{2\pi\i (x^B+x')k^B/N}  \cdot
e^{2\pi\i x'k'/N} \cdot
e^{2\pi\i x^B(k^B+k')/N}.
\]
Noticing that the first and the third terms depend only on $x'$ and
$k'$, respectively, we can factorize $M$ as $M=M_1 \cdot M_2 \cdot
M_3$, where $M_1$ and $M_3$ are diagonal matrices and $M_2$ is given
by
\begin{equation}
(M_2)_{x'k'} = e^{2\pi\i x'k'/N}
\label{eq:fracfft}
\end{equation}
for $0\le x',k' <s$. In fact, \eqref{eq:fracfft} is the matrix of the
fractional Fourier transform \cite{bailey-1991-ffta}, which can be
evaluated in only $O(s \log s)$ operations. Furthermore, since both
$M_1$ and $M_3$ are diagonal matrices, \eqref{eq:compB} can be
computed in $O(s \log s)$ steps as well.

Based on this observation, our algorithm takes the following form:
\begin{enumerate}
\item Construct a decomposition for $D = \{(x,k), k<c_x\}$.  Starting
  from $[0,N]^2$, we partition the boxes recursively. A box fully
  inside the $D$ is not further subdivided. The union of the boxes in
  the final decomposition is equal to $D$.
\item For $s=1,2,4,8,\cdots,N$, visit all of the boxes of size $s$ in
  the decomposition. Suppose $B$ is one such box. Compute the
  summation associated with $B$
  \[
  \sum_{k^B\le k < k^B+s} e^{2\pi\i xk/N} f_k
  \]
  for $x^B \le x <x^B+s$ using the fractional Fourier transform, and
  add the result to $\{u_x, x^B \le x <x^B+s\}$.
\end{enumerate}

The first step of our algorithm clearly takes at most $O(N\log N)$
steps. To estimate the complexity for the second step, one needs to have a bound on
the number of boxes of size $s$. Based on the construction of the
decomposition, we know that the center of a box $B$ of size $s$ is
at most of distance $s$ away from the curve $\{(x,c_x)\}$ because
otherwise $B$ would have been partitioned further. As a result, the
centers of all of the boxes of size $s$ must fall within a band
along $\{(x,c_x)\}$ of width $2s$. Noticing that $c^0(t)$ is smooth,
the length of $\{(x,c_x)\}$ is at most $N$. Therefore, the area of the
band is at most $O(Ns)$ and there are at most $O(Ns/s^2) = O(N/s)$
boxes of size $s$. Since we spend $O(s\log s)$ operations in the
fractional Fourier transform for each box of size $s$, the number of
steps for a fixed $s$ is $O(N/s \cdot s\log s) = O(N \log s) = O(N\log
N)$. Finally, summing over all $\log N$ possible values of $s$, we
conclude that our algorithm is $O(N \log^2 N)$. As no approximation
has been made, our algorithm is exact.


\begin{table}[ht]
  \begin{center}
    \includegraphics[height=2in]{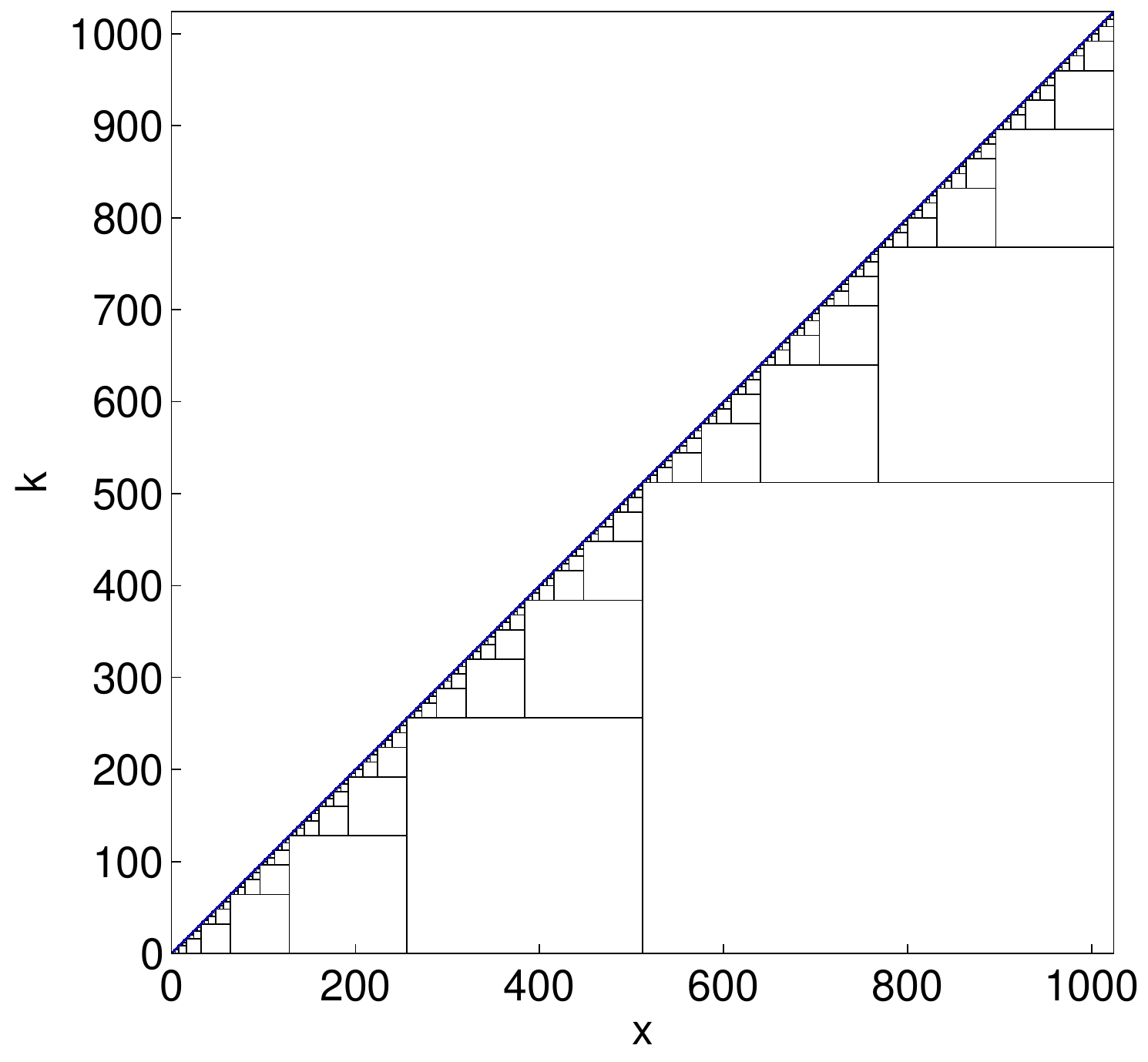}
    \includegraphics[height=2in]{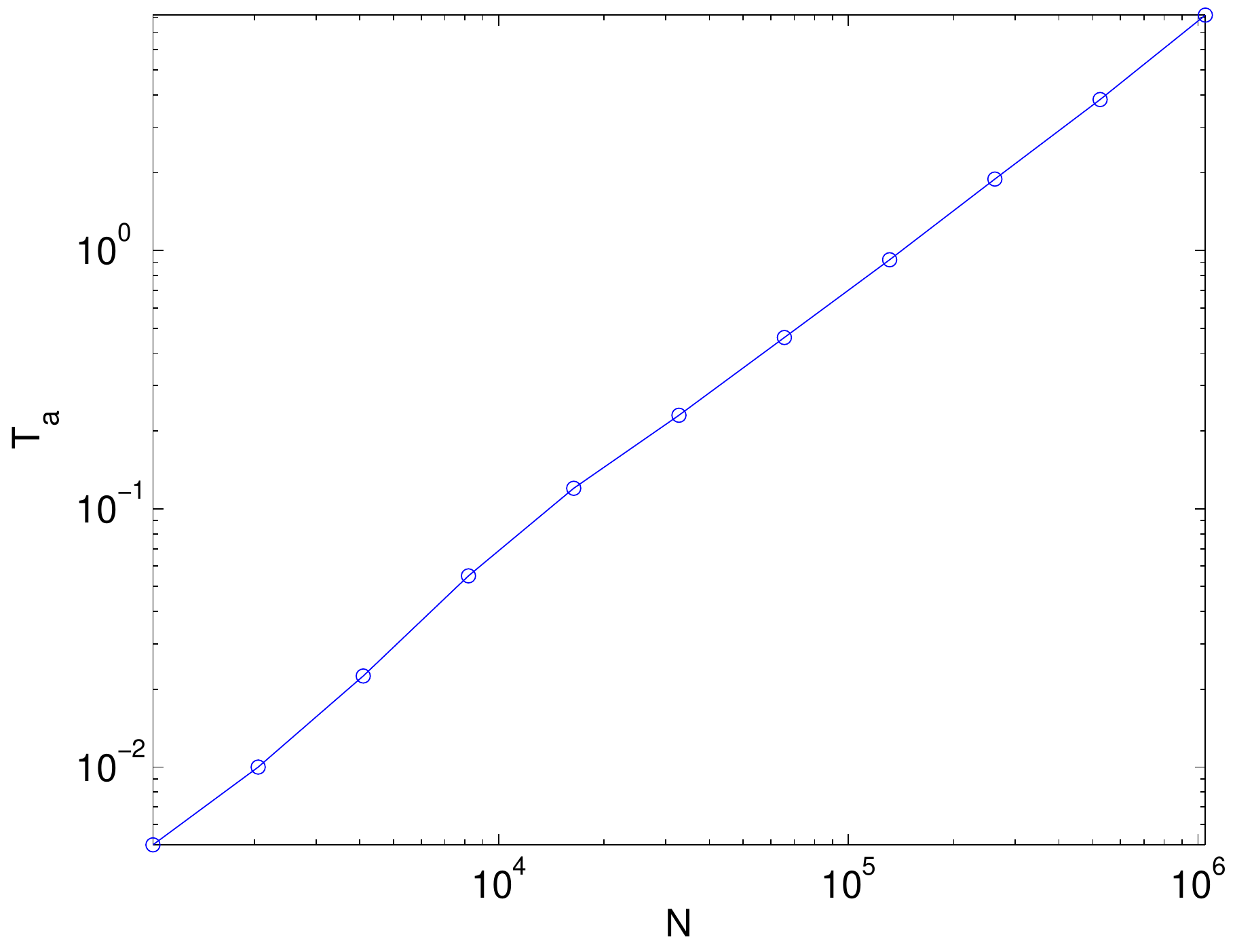}\\
    \vspace{0.1in}
    \begin{tabular}{|cccc|}
      \hline
      $N$ & $T_a$(sec) & $R_{d/a}$ & $R_{a/f}$ \\
      \hline
      1024 & 5.00e-03 & 3.00e+01 & 8.57e+00\\
      2048 & 1.00e-02 & 6.40e+01 & 9.27e+00\\
      4096 & 2.25e-02 & 1.14e+02 & 1.09e+01\\
      8192 & 5.50e-02 & 1.86e+02 & 1.36e+01\\
      16384 & 1.20e-01 & 3.41e+02 & 9.23e+00\\
      32768 & 2.30e-01 & 7.12e+02 & 4.15e+01\\
      65536 & 4.60e-01 & 1.42e+03 & 3.98e+01\\
      131072 & 9.20e-01 & 2.49e+03 & 3.68e+01\\
      262144 & 1.89e+00 & 4.16e+03 & 3.48e+01\\
      524288 & 3.84e+00 & 8.19e+03 & 8.30e+01\\
      1048576 & 8.16e+00 & 1.29e+04 & 7.25e+01\\
      \hline
    \end{tabular}
  \end{center}
  \caption{
    Top-left: the curve $c_x$ and the decomposition of $D$.
    Top-right: running time of our algorithm as a function of $N$.
    Bottom: the results for different values of $N$.
    $T_a$: the running time of our algorithm.
    $R_{d/a}$: the ratio between direct evaluation and our algorithm.
    $R_{a/f}$: the ratio between our algorithm and one execution of FFT of size $N$.
  }
  \label{tbl:1dlinear}
\end{table}

\begin{table}[ht]
  \begin{center}
    \includegraphics[height=2in]{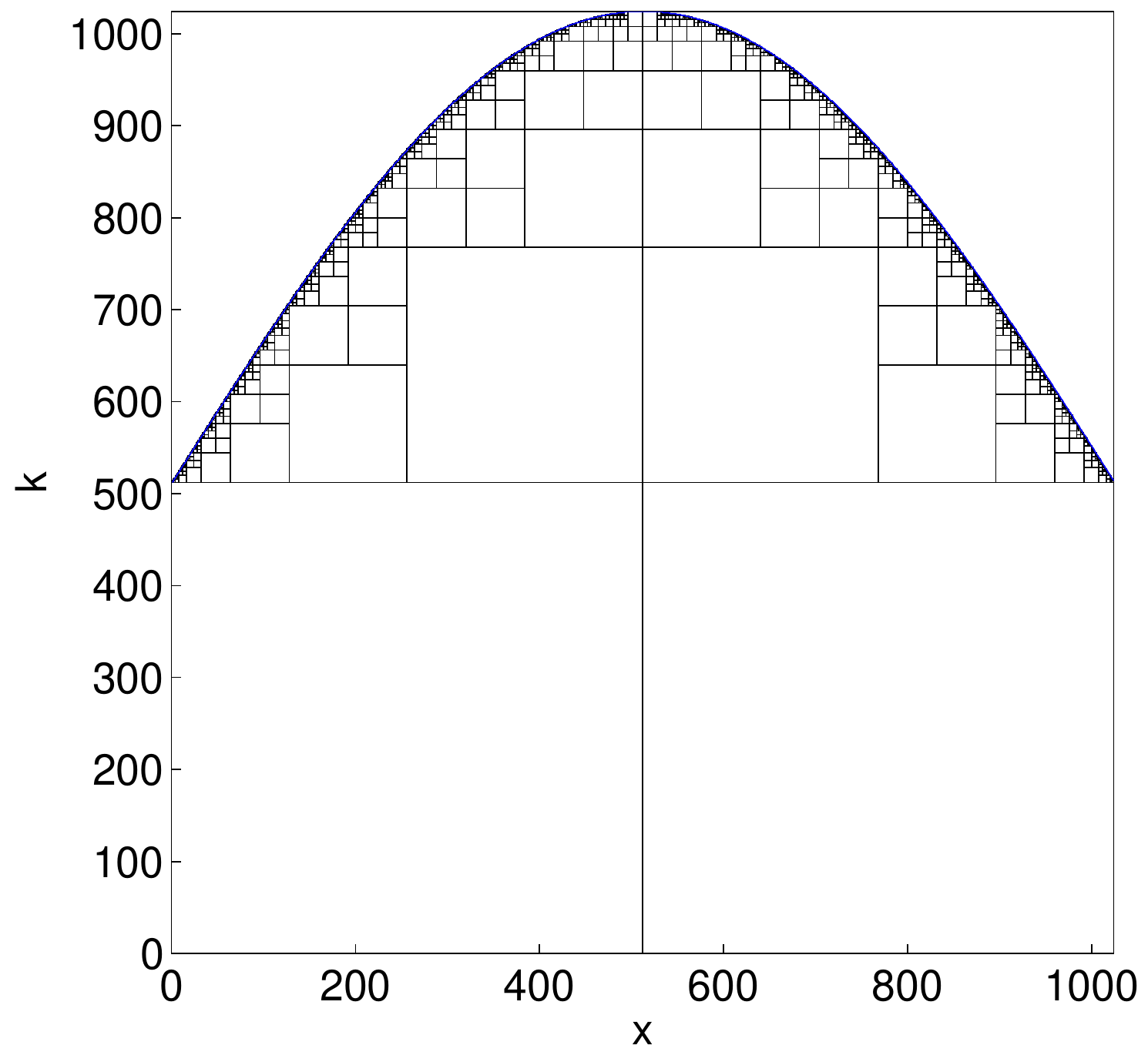}
    \includegraphics[height=2in]{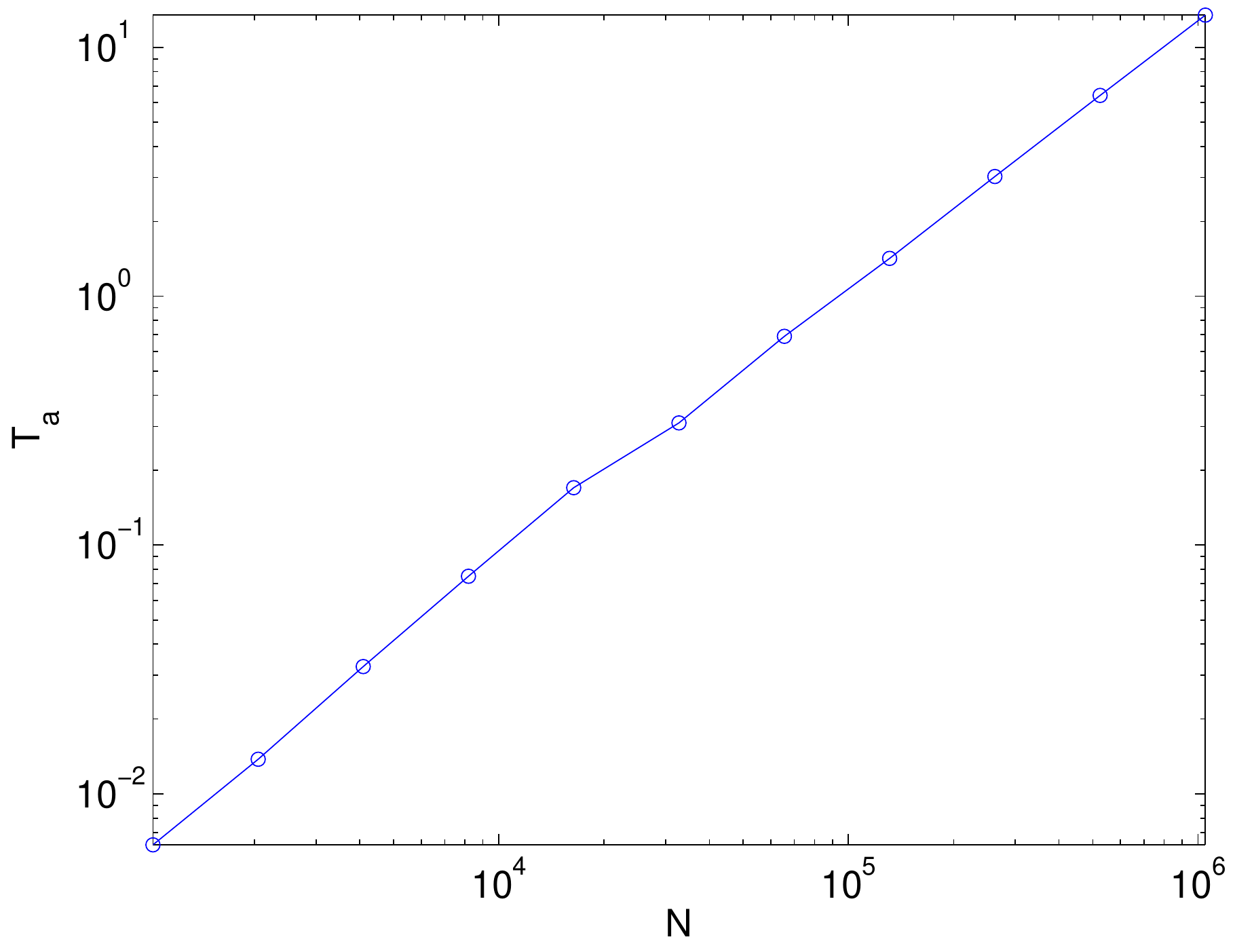}\\
    \vspace{0.1in}
    \begin{tabular}{|cccc|}
      \hline
      $N$ & $T_a$(sec) & $R_{d/a}$ & $R_{a/f}$ \\
      \hline
      1024 & 6.25e-03 & 2.40e+01 & 1.07e+01\\
      2048 & 1.38e-02 & 4.36e+01 & 1.27e+01\\
      4096 & 3.25e-02 & 7.88e+01 & 1.54e+01\\
      8192 & 7.50e-02 & 1.37e+02 & 1.86e+01\\
      16384 & 1.70e-01 & 2.41e+02 & 1.09e+02\\
      32768 & 3.10e-01 & 5.29e+02 & 9.92e+01\\
      65536 & 6.90e-01 & 8.90e+02 & 5.89e+01\\
      131072 & 1.42e+00 & 1.85e+03 & 5.83e+01\\
      262144 & 3.03e+00 & 2.81e+03 & 5.64e+01\\
      524288 & 6.42e+00 & 4.90e+03 & 1.17e+02\\
      1048576 & 1.35e+01 & 7.77e+03 & 1.28e+02\\
      \hline
    \end{tabular}
  \end{center}
  \caption{
    Top-left: the curve $c_x$ and the decomposition of $D$.
    Top-right: running time of our algorithm as a function of $N$.
    Bottom: the results for different values of $N$.
  }
  \label{tbl:1dsine}
\end{table}

\begin{table}[ht]
  \begin{center}
    \includegraphics[height=2in]{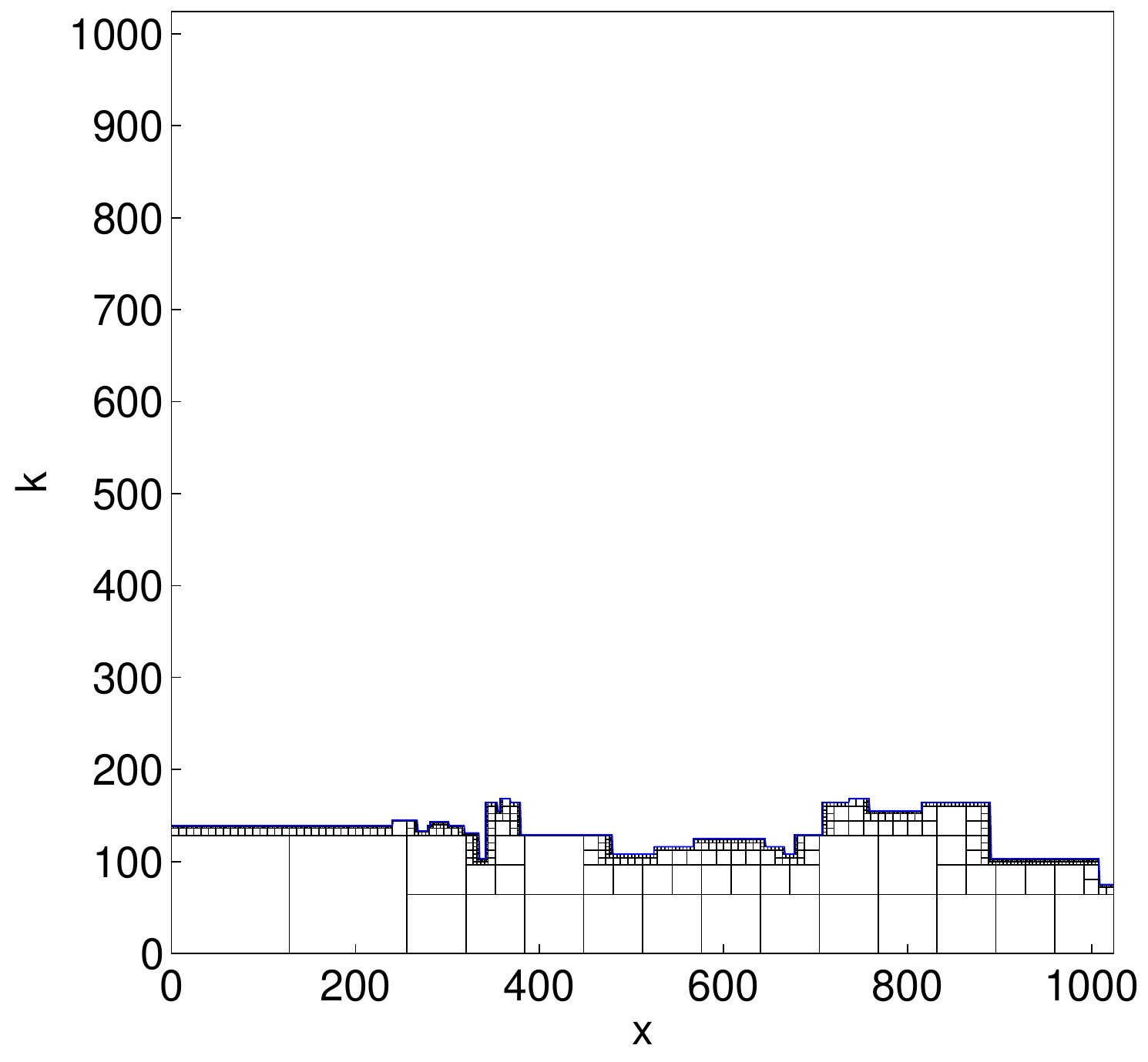}
    \includegraphics[height=2in]{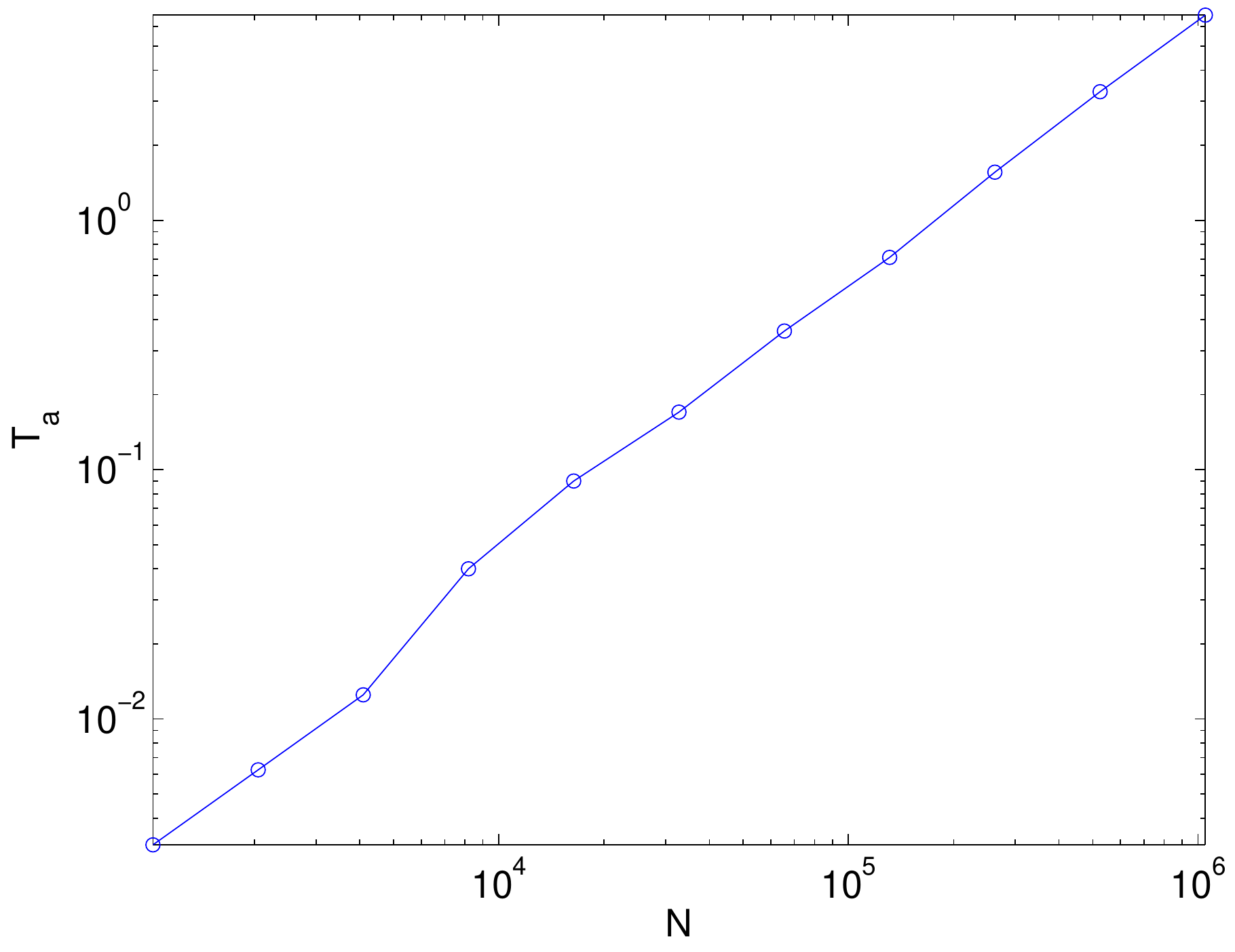}\\
    \vspace{0.1in}
    \begin{tabular}{|cccc|}
      \hline
      $N$ & $T_a$(sec) & $R_{d/a}$ & $R_{a/f}$ \\
      \hline
      1024 & 3.13e-03 & 4.48e+01 & 5.45e+01\\
      2048 & 6.25e-03 & 9.60e+01 & 6.56e+01\\
      4096 & 1.25e-02 & 1.92e+02 & 5.95e+00\\
      8192 & 4.00e-02 & 2.56e+02 & 5.00e+01\\
      16384 & 9.00e-02 & 4.27e+02 & 5.43e+00\\
      32768 & 1.70e-01 & 9.04e+02 & 9.85e+00\\
      65536 & 3.60e-01 & 1.82e+03 & 3.09e+01\\
      131072 & 7.10e-01 & 3.46e+03 & 2.91e+01\\
      262144 & 1.56e+00 & 5.46e+03 & 2.90e+01\\
      524288 & 3.28e+00 & 1.04e+04 & 6.64e+01\\
      1048576 & 6.66e+00 & 1.57e+04 & 6.50e+01\\
      \hline
    \end{tabular}
  \end{center}
  \caption{
    Top-left: the curve $c_x$ and the decomposition of $D$.
    Top-right: running time of our algorithm as a function of $N$.
    Bottom: the results for different values of $N$.
  }
  \label{tbl:1dkmax}
\end{table}


\subsection{Numerical results}

We apply our algorithm to several test examples to illustrate its
properties. All of the results presented here are obtained on a
desktop computer with a 2.8GHz CPU. For each example, we use the
following notations. $T_a$ is the running time of our algorithm in
seconds, $R_{d/a}$ is the ratio of the running time of direct
evaluation to $T_a$, and $R_{a/f}$ is the ratio of $T_a$ over the
running time of a Fourier transform (timed using FFTW
\cite{frigo-2005-difftw3}). As our algorithm is $O(N \log^2 N)$, we
expect $R_{d/a}$ to grow almost linearly and $R_{a/f}$ like $\log N$.

Tables \ref{tbl:1dlinear}, \ref{tbl:1dsine} and \ref{tbl:1dkmax}
summarize the results for three testing cases. The function in Table
\ref{tbl:1dkmax} corresponds to a 100 Hertz wave propagation through a
slice of the Marmousi velocity model \cite{versteeg-1994-me} taken at
2 km depth. From these examples, we observe clearly that $R_{d/a}$,
the ratio between the running times of direct evaluation and our
algorithm, indeed grows almost linearly in terms of $N$. Although the
ratio $R_{a/f}$ has some fluctuations, its value grows very slowly
with respect to $N$.

\section{Partial Fourier Transform in 2D}
\label{sec:2d}

A direct extension of the 1D algorithm to the 2D case would partition
the four dimensional summation domain $D = \{ (x,k), |k|<c_x\}$ with a
similar 4D tree structure. However, this does not result in an algorithm
with optimal complexity. To see this, let us count the number of boxes
of size $s$ in our tree structure. Repeating the argument used in
the complexity analysis of the 1D algorithm, we conclude that there
are about $N^3 s/s^4 = N^3/s^3$ boxes of size $s$. Even though the
computation associated with each box can be done in about $s^2 \log s$
steps, the total operation count for a fixed $s$ is about $N^3/s^3
\cdot s^2 \log s = N^3/s \log s$, which is much larger than the degree
of freedom $N^2$ for small values of $s$.

\subsection{Algorithm description}

Noticing that only $|k|$ appears in the constraint of the 2D partial
Fourier transform
\[
u_x = \sum_{|k|<c_x} e^{2\pi\i x\cdot k/N} f_k \quad
0 \le x_1,x_2 < N,
\]
we study a different set $R = \{ (x,r), r<c_x\}$ instead.

\begin{figure}[ht]
  \begin{center}
    \includegraphics[height=2in]{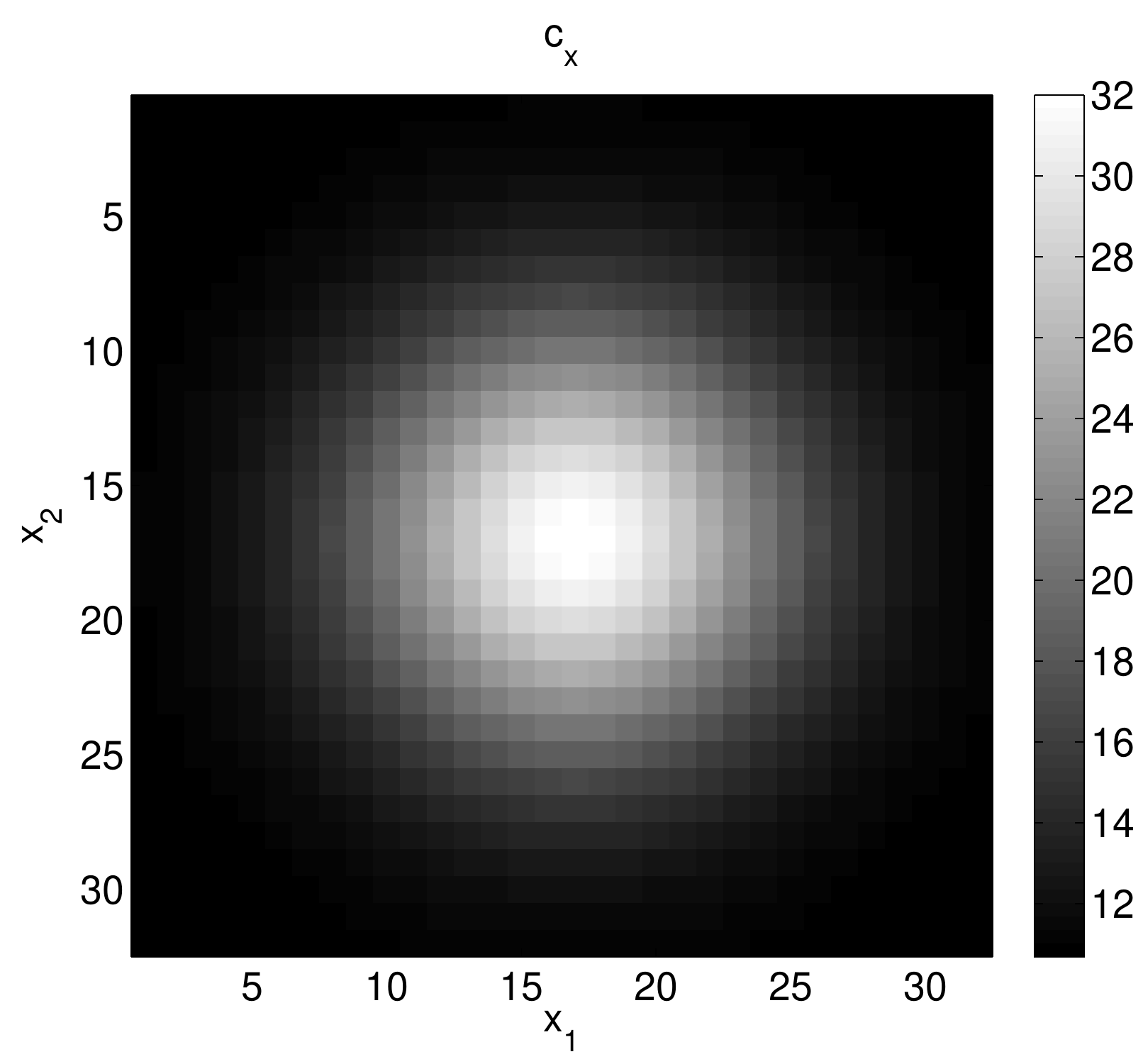}
    \includegraphics[height=2in]{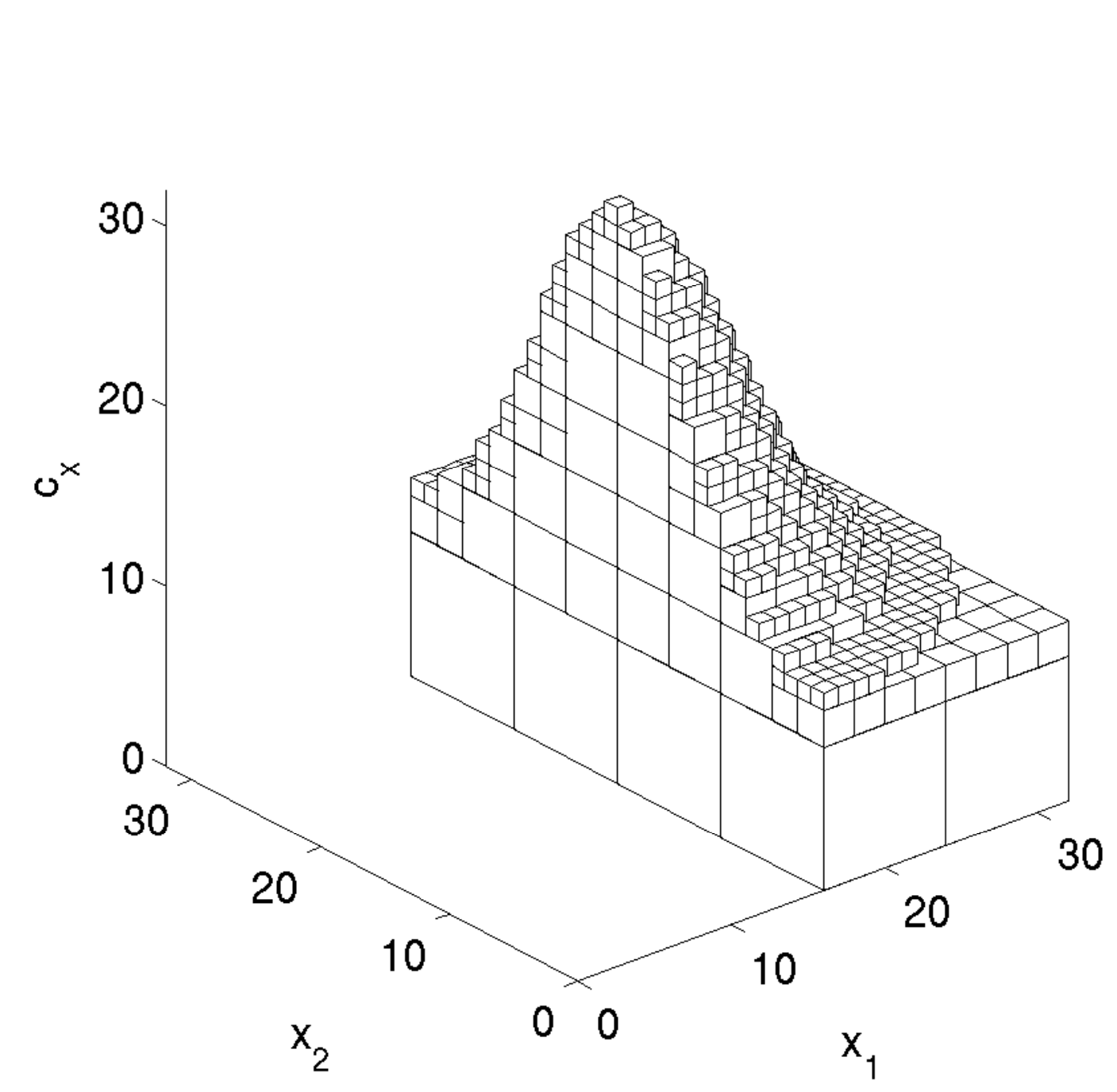}
  \end{center}
  \caption{ Left: $c_x$ is a Gaussian function. Right: a cross section
    view of the multiscale decomposition of $R$ (the domain below the
    surface $c_x$.  Here $N=32$.  }
  \label{fig:decomp2d}
\end{figure}

The algorithm first generates a decomposition for $R$. Similar to the
1D case, we partition the box $[0,N]^3$ through recursive subdivision.
A box is not further subdivided if it fully resides in $R$. The union
of all the boxes inside the decomposition is exactly the set $R$
(see Figure \ref{fig:decomp2d}).

\begin{figure}[ht]
  \begin{center}
    \includegraphics[height=1.8in]{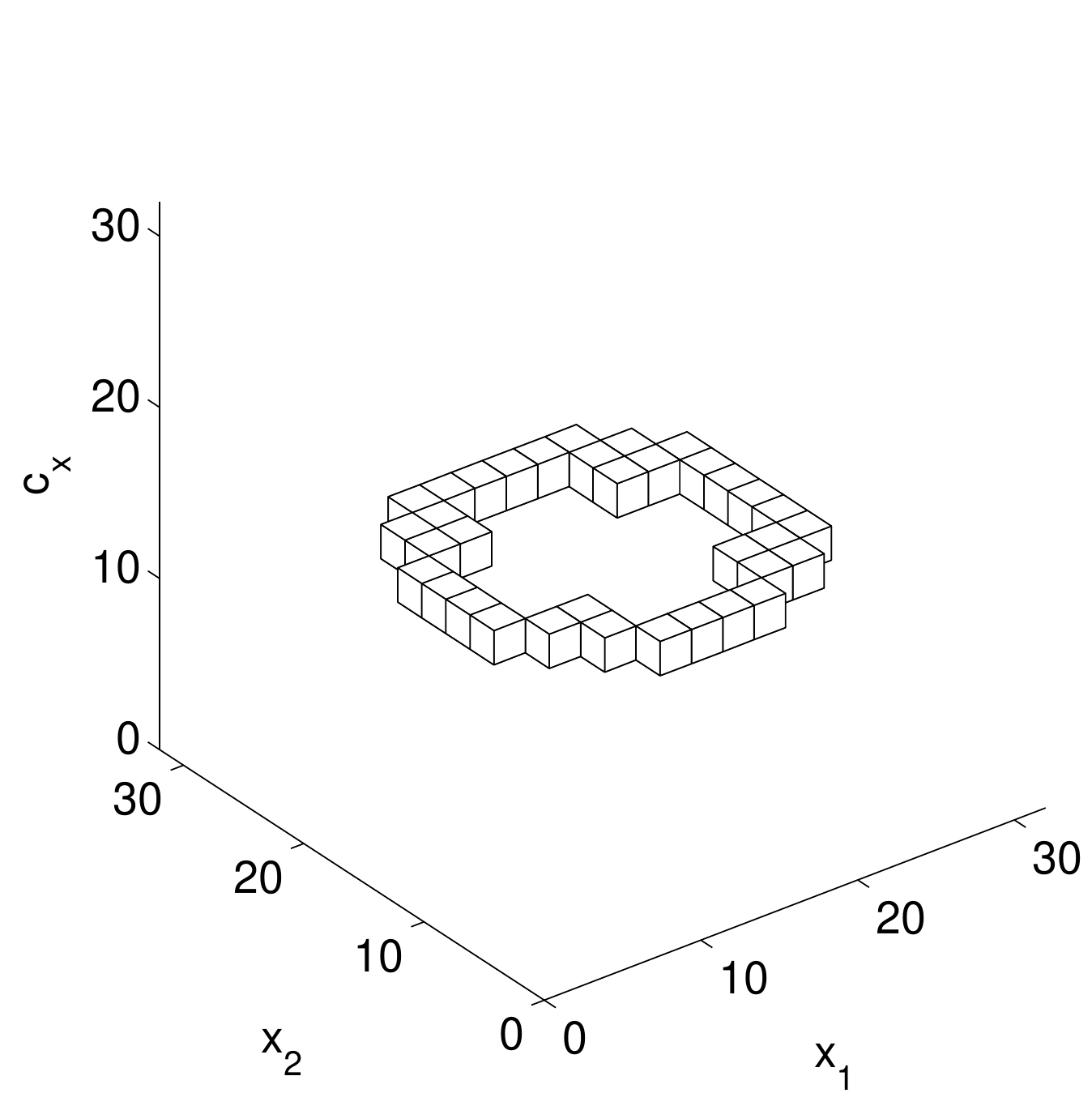}
    \includegraphics[height=1.8in]{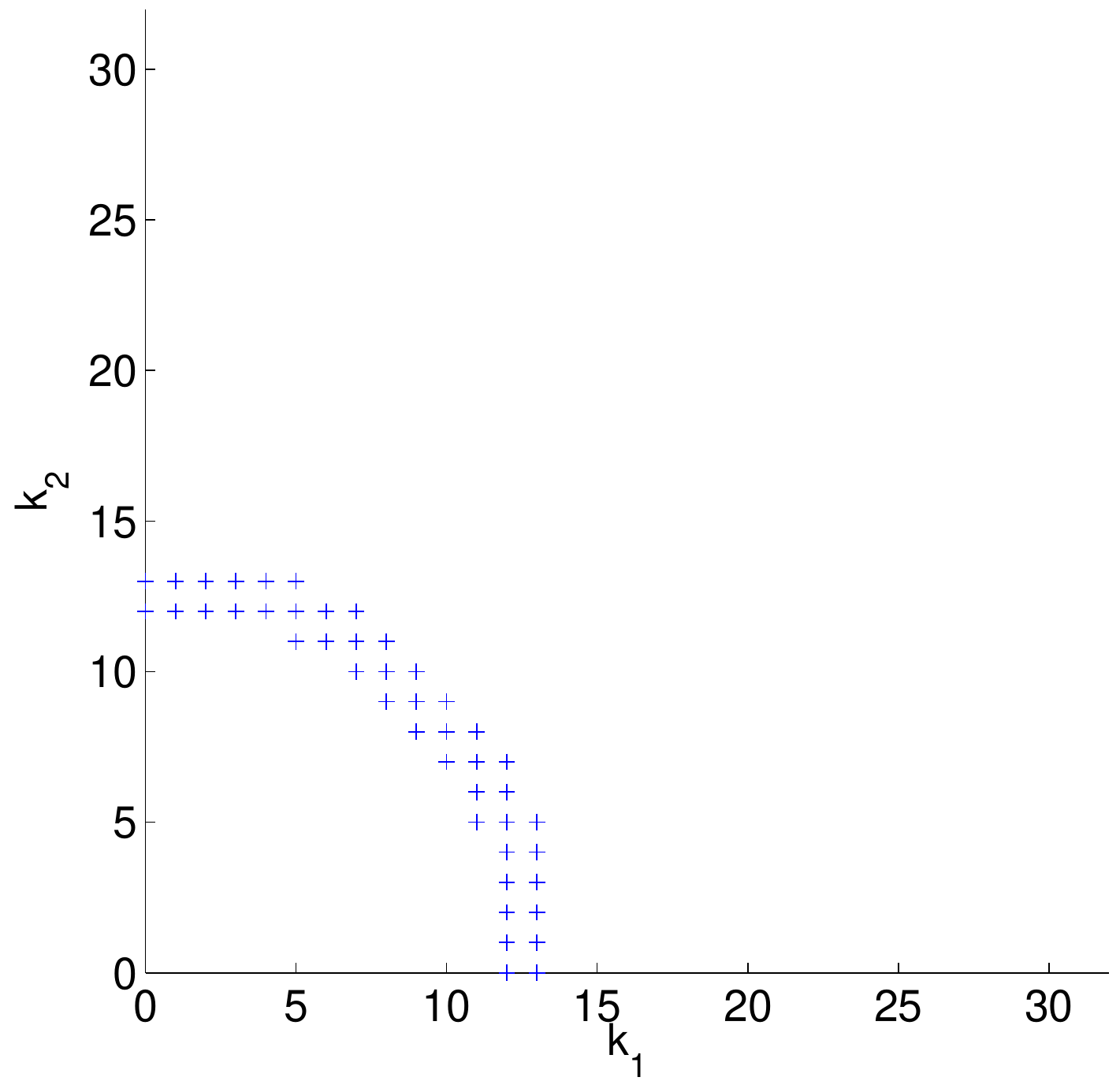}
    \includegraphics[height=1.8in]{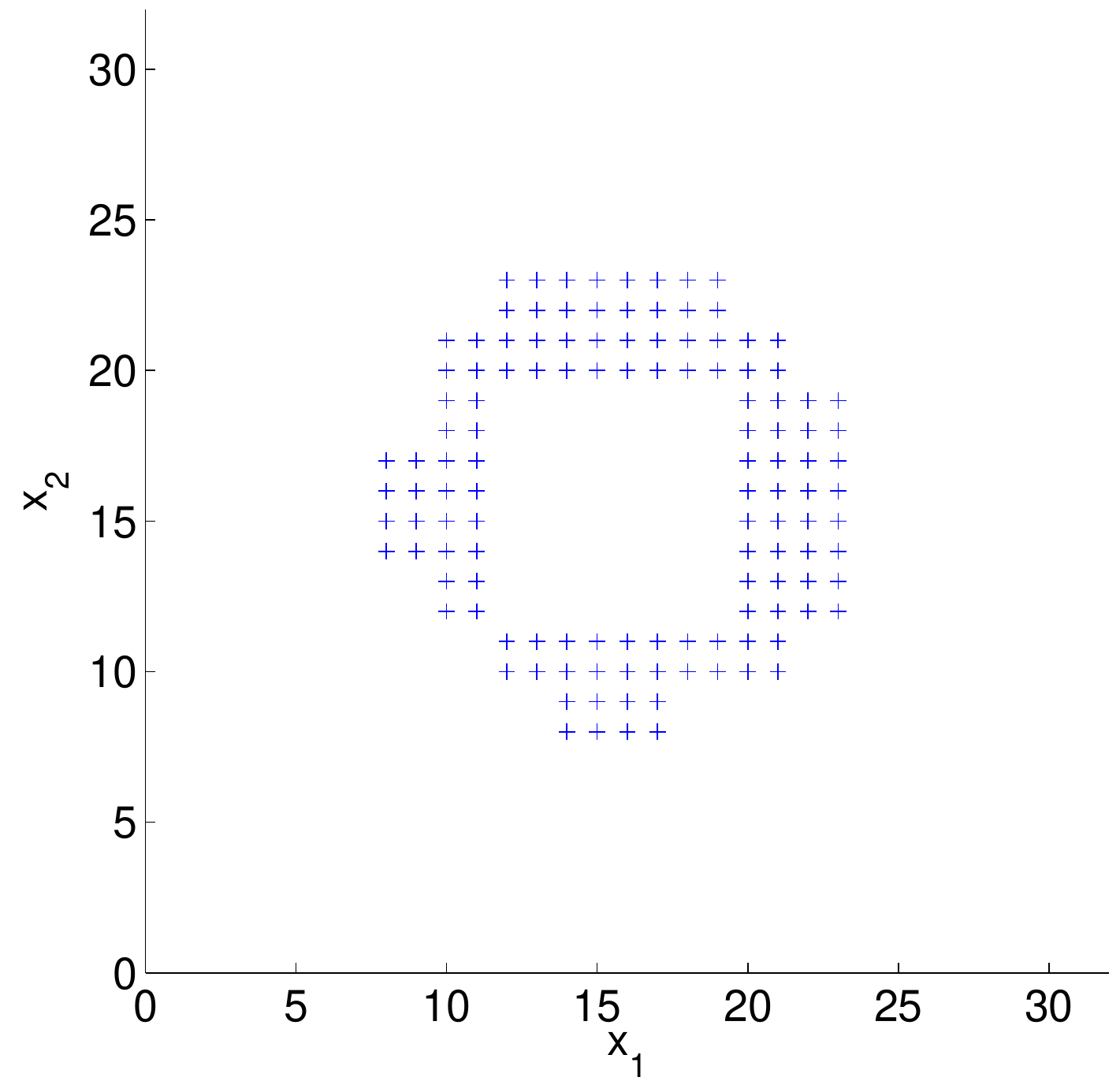}
  \end{center}
  \caption{ Left: the boxes inside the set $G^A$ for an given interval
    $A$ of the $r$ coordinate.  Middle: the set $K^A$ in the
    $(k_1,k_2)$ plane.  Right: the set $X^A$ in the $(x_1,x_2)$ plane.
    The butterfly algorithm in \cite{ying-2008-sftba} evaluates
    the summation between $X^A$ and $K^A$ in almost linear time.
    Here $N=32$.
  }
  \label{fig:XAKA}
\end{figure}

The projection of any box of our decomposition onto the $r$ coordinate
is a dyadic interval. Let us consider one such interval $A$ of size
$s$ and denote $G^A$ to be the set of all cubes that project onto $A$.
We define $K^A$ to be the set $\{k, |k|\in A\}$ and $X^A$ to be the
image of the points in $G^A$ under the projection onto the $(x_1,x_2)$
plane. Since $A$ is an interval of size $s$, $K^A$ is in fact a band
in the $(k_1,k_2)$ domain with length $O(N)$ and width $s$. Noticing
that the surface $c^0(t)$ used to define $\{c_x,0\le x_1,x_2 <N\}$ is
smooth, the set $X_A$ is also a band in the $(x_1,x_2)$ domain with
length $O(N)$ and width $O(s)$ (see Figure \ref{fig:XAKA}).

The part of the summation associated with the interval $A$ is
\begin{equation}
  \sum_{k\in K^A} e^{2\pi\i x\cdot k /N} f_k
  \label{eq:compA}
\end{equation}
for $x\in X^A$. Since $X^A$ and $K^A$ are two bands in $[0,N)^2$,
\eqref{eq:compA} is indeed a Fourier transform problem with sparse
data. To compute \eqref{eq:compA}, we utilize the solution proposed in
\cite{ying-2008-sftba}. This approach is a butterfly algorithm based
on \cite{michielssen-1996-mmda,oneil-2007-ncabft} and computes an
approximation of \eqref{eq:compA} in $O(\max(|X^A|,|K^A|)
\log(\max(|X^A|,|K^A|))$ operations, almost linear in terms of the
degree of freedom.

Combining these ideas, we have the following algorithm:
\begin{enumerate}
\item Construct a decomposition for $R = \{(x,r), r<c_x\}$.  Starting
  from $[0,N]^3$, we partition the boxes recursively. A box fully
  inside the $R$ is not further subdivided. The union of the boxes in
  the final decomposition is equal to $R$.
\item For $s=1,2,4,8,\cdots,N$, visit all the dyadic intervals of size
  $s$ in the $r$ coordinate. Suppose that $A$ is one such interval.
  Compute the summation associated with $A$
  \[
  \sum_{k\in K^A} e^{2\pi\i x\cdot k /N} f_k
  \]
  for $x\in X^A$ using the butterfly procedure in
  \cite{ying-2008-sftba}, and add the result to $\{u_x, x\in X^A\}$.
\end{enumerate}

Let us consider now the complexity of this algorithm. The first step
of our algorithm takes only $O(N^2\log N)$ steps. In order to estimate
the number of operations used in the second step, let us consider a
fixed $s$. For each interval $A$ of size $s$, the number of steps used
in $O(\max(|X^A|,|K^A|) \log(\max(|X^A|,|K^A|)) = O(|X^A|\log |X^A|) +
O(|K^A|\log |K^A|)$. Summing over all boxes of size $s$, we get
\[
O\left( \sum_{A: |A|=s} |X^A|\log |X^A| \right) + 
O\left( \sum_{A: |A|=s} |K^A|\log |K^A| \right).
\]
Noticing $\sum_{A: |A|=s} |X^A| = \sum_{A: |A|=s} |K^A| = N^2$, the
above quantity is clearly bounded by $O(N^2\log N)$. Finally, after
summing over all different values of $s$, we have a total complexity
of order $O(N^2 \log^2 N)$.

\subsection{Numerical results}

We apply our algorithm to several examples in this section. In
\cite{ying-2008-sftba}, equivalent charges located at Cartesian grids
are used as the low rank representations in the butterfly algorithm to
control the accuracy of the method. The size of the Cartesian grid $p$
controls the accuracy of our algorithm. Here, we pick $p$ to be 5 or
9. To quantify the error, we select a set $S \subset
\{0,1,\cdots,N-1\}^2$ of size 100 and estimate the error by
\[
\sqrt{
  \frac{ \sum_{x\in S} |u_x - u_x^a |^2 } { \sum_{x\in S} |u_x|^2 }
}
\]
where $\{u_x\}$ are the exact results and $\{u_x^a\}$ are our
approximations.

\begin{table}[ht]
  \begin{center}
    \includegraphics[height=2in]{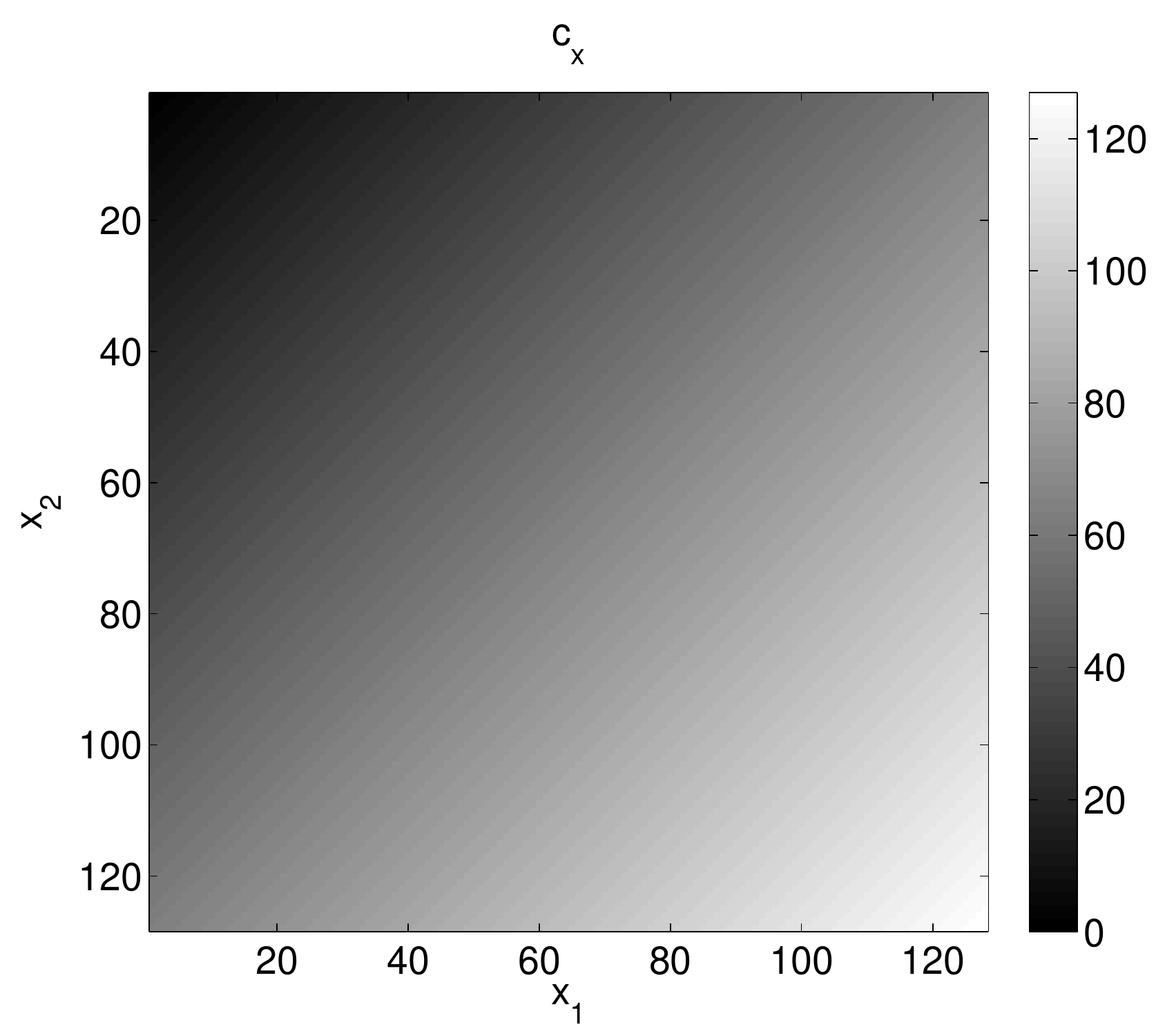}
    \includegraphics[height=2in]{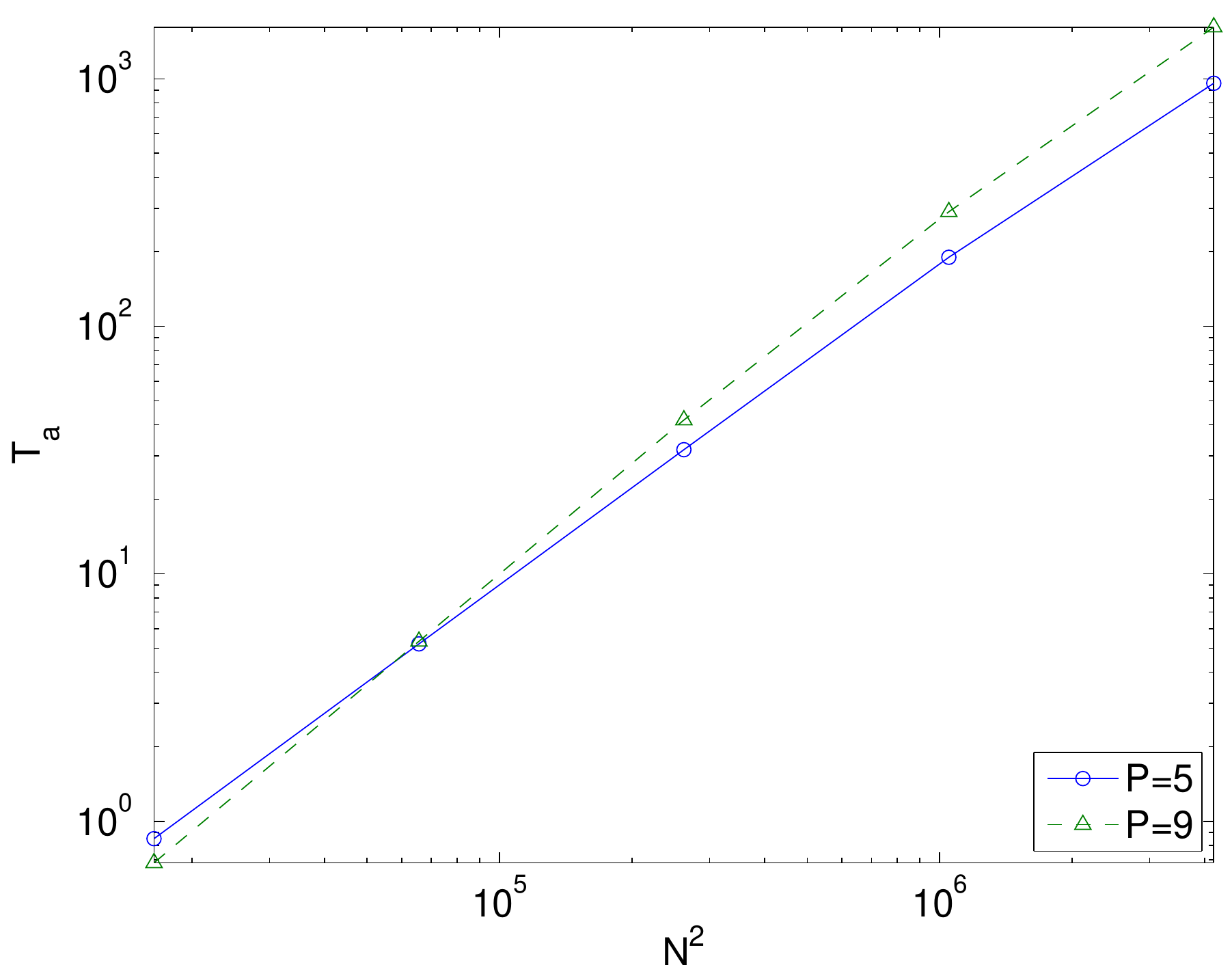}\\
    \vspace{0.1in}
    \begin{tabular}{|ccccc|}
      \hline
      $(N,p)$ & $T_a$(sec) & $R_{d/a}$ & $R_{a/f}$ & $E_a$ \\
      \hline
      (128,5) & 8.50e-01 & 4.82e+01 & 9.89e+02 & 5.25e-04\\
      (256,5) & 5.21e+00 & 1.26e+02 & 1.59e+03 & 9.85e-04\\
      (512,5) & 3.17e+01 & 3.76e+02 & 1.56e+03 & 7.40e-04\\
      (1024,5)& 1.90e+02 & 1.13e+03 & 1.61e+03 & 1.21e-03\\
      (2048,5)& 9.60e+02 & 4.94e+03 & 1.31e+03 & 1.20e-03\\
      \hline
      (128,9) & 6.80e-01 & 6.02e+01 & 7.25e+02 & 1.46e-14\\
      (256,9) & 5.33e+00 & 1.29e+02 & 1.42e+03 & 8.87e-09\\
      (512,9) & 4.18e+01 & 2.98e+02 & 1.63e+03 & 2.17e-08\\
      (1024,9)& 2.90e+02 & 7.49e+02 & 1.97e+03 & 9.34e-09\\
      (2048,9)& 1.62e+03 & 2.92e+03 & 2.10e+03 & 2.39e-08\\
      \hline
    \end{tabular}
  \end{center}
  \caption{
    Top-left: $c_x$ when $N=128$.
    Top-right: running time of our algorithm as a function of $N^2$.
    Bottom: the results for different values of $N$.
    $T_a$: the running time of our algorithm.
    $R_{d/a}$: the ratio between direct evaluation and our algorithm.
    $R_{a/f}$: the ratio between our algorithm and one execution of FFT of size $N$.
    $E_a$: the estimated error.
  }
  \label{tbl:2dplne}
\end{table}

\begin{table}[ht]
  \begin{center}
    \includegraphics[height=2in]{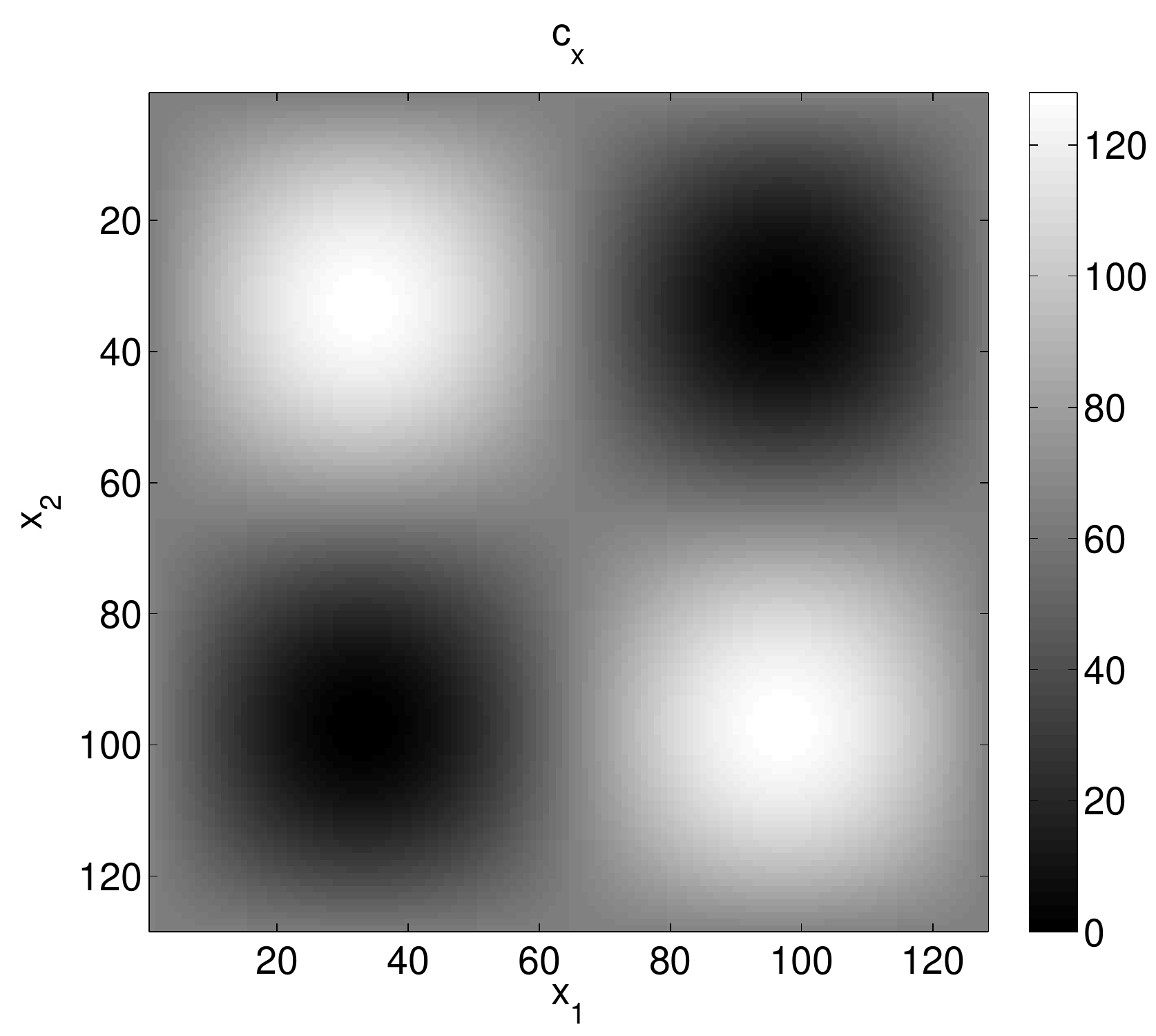}
    \includegraphics[height=2in]{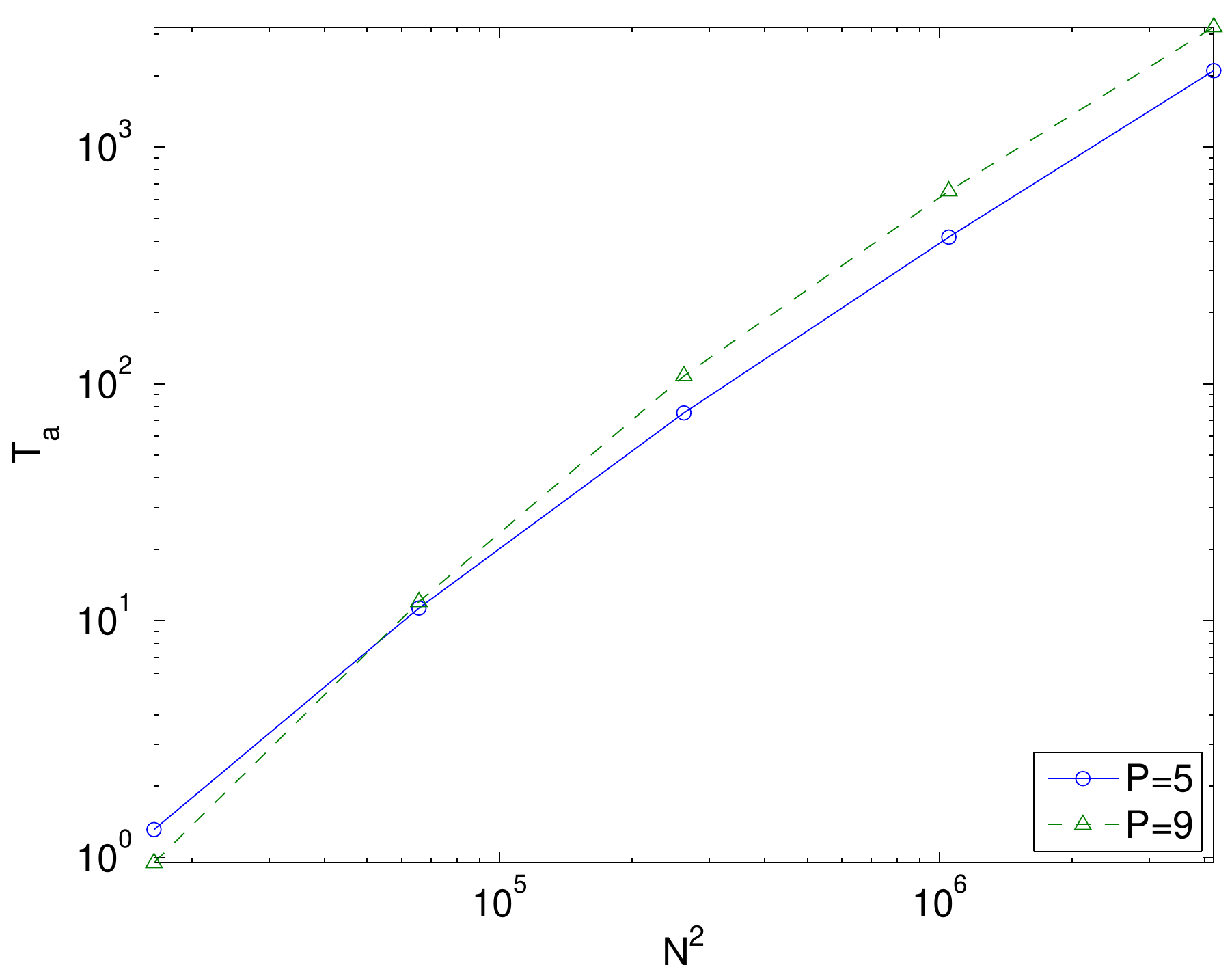}\\
    \vspace{0.1in}
    \begin{tabular}{|ccccc|}
      \hline
      $(N,p)$ & $T_a$(sec) & $R_{d/a}$ & $R_{a/f}$ & $E_a$ \\
      \hline
      (128,5) & 1.31e+00 & 3.13e+01 & 1.52e+03 & 7.11e-04\\
      (256,5) & 1.13e+01 & 6.35e+01 & 2.91e+03 & 5.17e-04\\
      (512,5) & 7.53e+01 & 1.65e+02 & 2.94e+03 & 1.27e-03\\
      (1024,5)& 4.16e+02 & 5.45e+02 & 2.62e+03 & 8.32e-04\\
      (2048,5)& 2.10e+03 & 2.28e+03 & 2.71e+03 & 8.88e-04\\
      \hline
      (128,9) & 9.50e-01 & 4.31e+01 & 1.01e+03 & 2.02e-14\\
      (256,9) & 1.20e+01 & 5.75e+01 & 3.06e+03 & 6.73e-09\\
      (512,9) & 1.08e+02 & 1.17e+02 & 4.20e+03 & 9.35e-09\\
      (1024,9)& 6.51e+02 & 3.46e+02 & 4.10e+03 & 1.32e-08\\
      (2048,9)& 3.21e+03 & 1.48e+03 & 4.12e+03 & 2.00e-08\\
      \hline
    \end{tabular}
  \end{center}
  \caption{
    Top-left: $c_x$ when $N=128$.
    Top-right: running time of our algorithm as a function of $N^2$.
    Bottom: the results for different values of $N$.
  }
  \label{tbl:2dsine}
\end{table}

\begin{table}[ht]
  \begin{center}
    \includegraphics[height=2in]{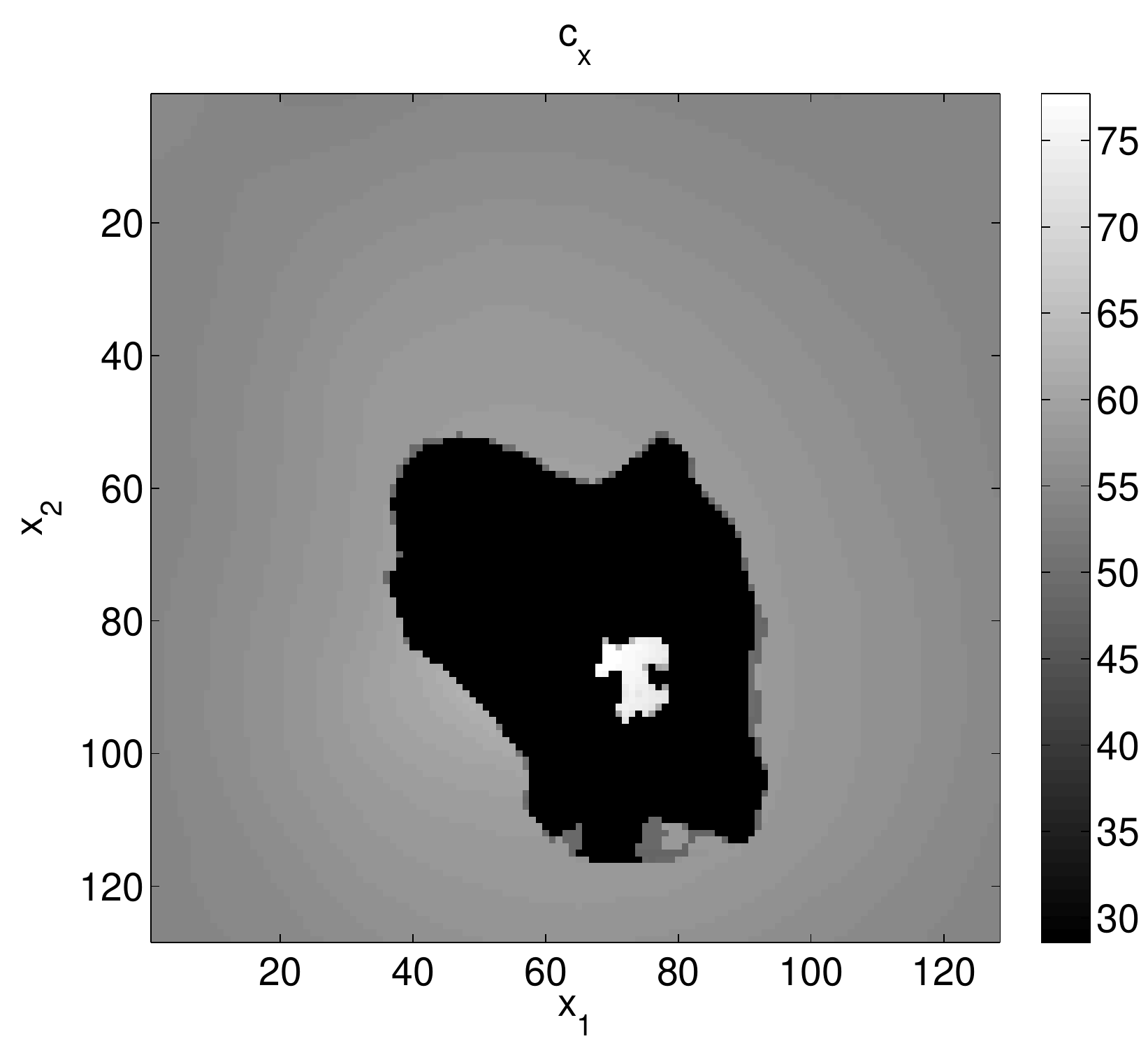}
    \includegraphics[height=2in]{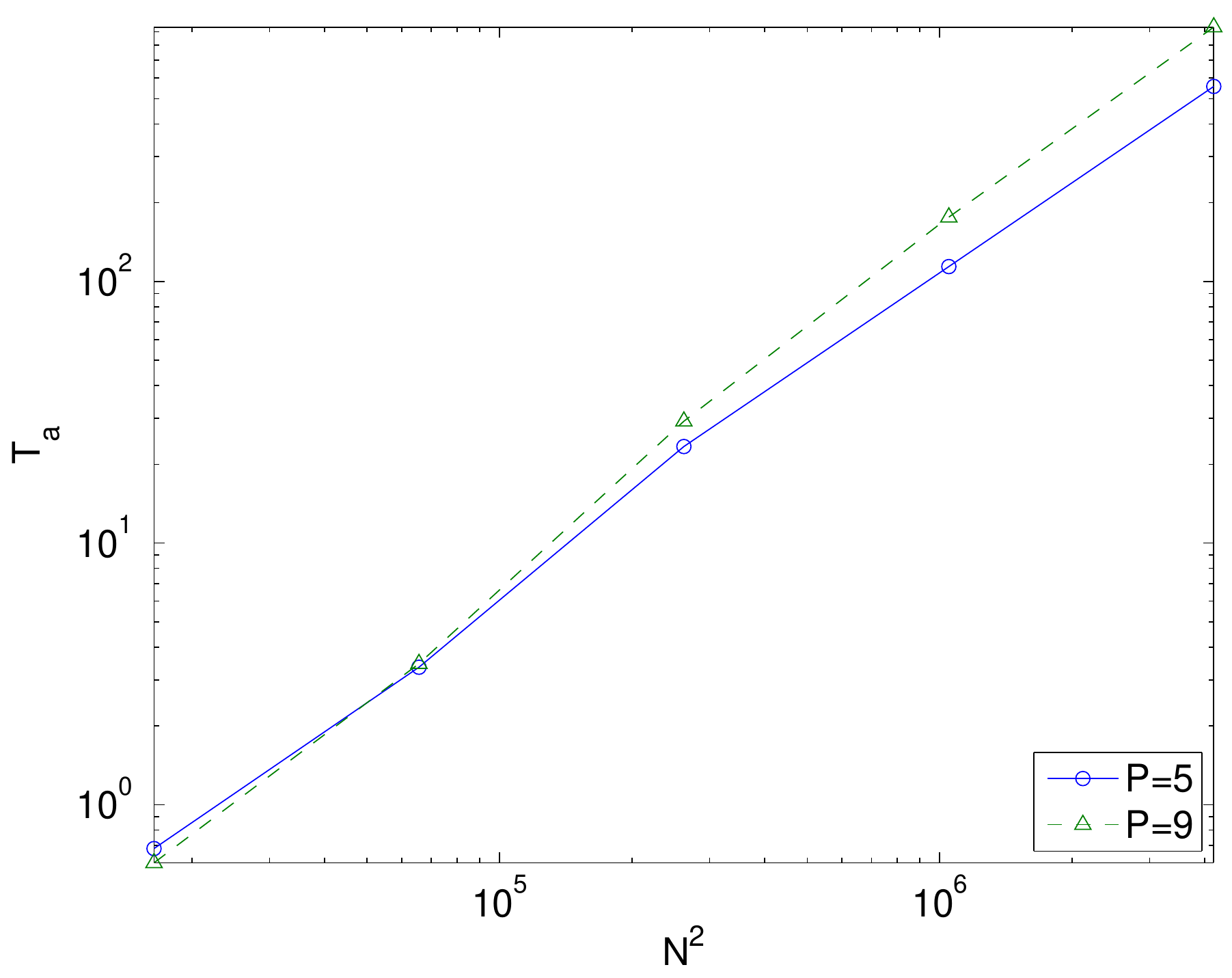}\\
    \vspace{0.1in}
    \begin{tabular}{|ccccc|}
      \hline
      $(N,p)$ & $T_a$(sec) & $R_{d/a}$ & $R_{a/f}$ & $E_a$ \\
      \hline
      (128,5) & 6.80e-01 & 6.02e+01 & 7.25e+02 & 4.81e-04\\
      (256,5) & 3.35e+00 & 2.15e+02 & 8.58e+02 & 7.17e-04\\
      (512,5) & 2.34e+01 & 5.72e+02 & 6.18e+02 & 9.29e-04\\
      (1024,5) & 1.14e+02 & 2.02e+03 & 6.79e+02 & 1.04e-03\\
      (2048,5) & 5.57e+02 & 8.83e+03 & 6.59e+02 & 1.07e-03\\
      \hline
      (128,9) & 6.00e-01 & 6.83e+01 & 6.40e+02 & 7.66e-15\\
      (256,9) & 3.46e+00 & 1.99e+02 & 8.52e+02 & 1.13e-08\\
      (512,9) & 2.92e+01 & 4.09e+02 & 1.24e+03 & 1.68e-08\\
      (1024,9) & 1.76e+02 & 1.34e+03 & 1.14e+03 & 2.07e-08\\
      (2048,9) & 9.39e+02 & 5.03e+03 & 1.21e+03 & 2.62e-08\\
      \hline
    \end{tabular}
  \end{center}
  \caption{
    Top-left: $c_x$ when $N=128$.
    Top-right: running time of our algorithm as a function of $N^2$.
    Bottom: the results for different values of $N$.
  }
  \label{tbl:2dkmax2}
\end{table}


Similar to the 1D case, the following notations are used: $T_a$ is the
running time of our algorithm in seconds, $R_{d/a}$ is the ratio of
the running time of direct evaluation to $T_a$, $R_{a/f}$ is the ratio
of $T_a$ over the running time of a Fourier transform (timed using
FFTW \cite{frigo-2005-difftw3}), and finally $E_a$ is the estimated
error.

The numerical results are summarized in Tables \ref{tbl:2dplne},
\ref{tbl:2dsine} and \ref{tbl:2dkmax2}. The function in Table
\ref{tbl:2dkmax2} corresponds to a 50 Hertz wave propagation through a
slice of the SEG/EAGE velocity model \cite{aminzadeh-1997-3som} taken
at 1.5 km depth. From these numbers, we see that our implementation
indeed has a complexity almost linear in terms of the number of grid
points. Due to the complex structure of the butterfly procedure, the
constant of our algorithm is quite large compared to the one of FFTW.

\section{Conclusions and Discussions} 
\label{sec:conc}

In this paper, we introduced two efficient algorithms for computing
partial Fourier transforms in one and two dimensions. In both cases,
we start by decomposing the appropriate summation domain in a
multiscale way into simple pieces and apply existing fast algorithms
on each piece to get optimal efficiency. In 1D, the fractional Fourier
transform is used. In 2D, we resort to the butterfly algorithm for
sparse Fourier transform proposed in \cite{ying-2008-sftba}. As a
result, both of our algorithms are asymptotically only $O(\log N)$
times more expensive than the FFT.

In Tables \ref{tbl:2dplne}, \ref{tbl:2dsine} and \ref{tbl:2dkmax2}, we
notice that our 2D algorithm has a relatively large constant. One
obvious direction of future research is to improve on our current
implementation of the butterfly algorithm. Another alternative is to
seek different ways for computing the Fourier transforms with sparse
data.  In the past several years, several algorithms have been
developed to address similar oscillatory behavior efficiently (see,
for example,
\cite{averbuch-2000-ecoi,cheng-2006-awfmm,engquist-2007-fdmaok}). It
would be interesting to see whether these ideas can be used in the
setting of the Fourier transform with sparse data.

As we mentioned earlier, this research is motivated mostly by the wave
extrapolation algorithm in seismic imaging. Our model problem
considers only one of the challenges, i.e., the existence of the
summation constraint. The other challenge is to improve the evaluation
of the $e^{2\pi\i \sqrt{ \omega^2/v^2(x) - k^2 }\cdot z}$ term, for
example, by approximating it on each of the simple summation
components. Research along this direction will be presented in a
future report.

{\bf Acknowledgments.} The first author is partially supported by an
Alfred P. Sloan Fellowship and the startup grant from the University
of Texas at Austin.

\bibliographystyle{abbrv}
\bibliography{ref}

\end{document}